\documentclass[11pt]{article}

\input{amssym.def}
\input{amssym.tex}

\usepackage{amsthm}
\usepackage{epsfig}
\usepackage{color}

\title{LERF and the Lubotzky-Sarnak Conjecture}
\author{M. Lackenby\thanks{First author supported by
an EPSRC Advanced Research Fellowship}, 
D. D. Long \& A. W. Reid\thanks{Last two authors supported in part by the NSF}}

\def\vol{{\rm vol}}
\def\inf{{\rm inf}}

\def\SO{{\rm SO}}
\def\Sp{{\rm Sp}}

\def\qed{ $\sqcup\!\!\!\!\sqcap$}

\def\Area{\mbox{{\rm Area}}}

\newtheorem{theorem}{Theorem}[section]

\newtheorem{corollary}[theorem]{Corollary}
\newtheorem{proposition}[theorem]{Proposition}
\newtheorem{conjecture}[theorem]{Conjecture}

\begin{document}
\maketitle

%
%
%
%
\section{Introduction}
We begin by recalling the definition of Property $\tau$.
Let $X$ be a finite graph,
and let $V(X)$ denote its vertex set. For any subset $A$ of
$V(X)$, let $\partial A$ denote those edges with one endpoint
in $A$ and one not in $A$. Define the {\em Cheeger constant}
of $X$ to be
$$h(X) = \min \left \{ {|\partial A| \over |A|} :
A \subset V(X) \hbox{ and } 0 < |A| \leq |V(X)|/2 \right \}.$$
Now let $G$ be a group with a finite symmetric generating set $S$.
For any subgroup $G_i$ of $G$, let $X(G/G_i;S)$ be the
Schreier coset graph of $G/G_i$ with respect to $S$.
Then $G$ is said to have {\em Property $\tau$ with
respect to a collection of finite index subgroups} $\{ G_i \}$
if $\inf_i \;h(X(G/G_i; S)) > 0$. This turns out not to depend on
the choice of finite generating set $S$.  Also, $G$ is said to have 
{\em Property 
$\tau$} if it has Property $\tau$ with respect to the collection of
all subgroups of finite index in $G$.

In the context of of finite volume hyperbolic manifolds,
Lubotzky and Sarnak made the following conjecture.

\begin{conjecture}
The fundamental group of any finite volume hyperbolic n-manifold
does not have Property $\tau$.\end{conjecture}

It is easy to check that if a group $G$ contains a finite
index subgroup surjecting onto $\bf Z$, then $G$ does not have Property
$\tau$, and it is this that has attracted attention to the Lubotzky-Sarnak
conjecture recently. This is particularly relevant in the context
of hyperbolic 3-manifolds (see \cite{La1} and
\cite{LaLR} for example), in part due to the connection to the
virtual positive first Betti number conjecture from 3-manifold
topology (see \cite{LaLR} for a  discussion). 

While it appears to be much weaker than the virtual positive first Betti
number conjecture, it appears that there is no method known to show that the
fundamental group of a finite volume hyperbolic $n$-manifold does not
have Property $\tau$ without exhibiting a surjection onto $\bf Z$ from a
finite index subgroup.  The main result of this note provides a
method for hyperbolic 3-manifolds. 

To state the main result we require some additional terminology.

Let $G$ be a finitely generated group and $H$ a finitely generated
subgroup.  $G$ is {\em $H$-separable} if $H$ is closed in the
profinite topology on $G$, and $G$ is called {\em LERF or subgroup
  separable} if $G$ is $H$-separable for every finitely generated
subgroup $H<G$. 
We say that $H$ is {\em engulfed in $G$} if there is a
proper finite index subgroup $K<G$ with $H<K$.  In the context of
hyperbolic $3$-manifolds, it turns out that these two notions are
intimately related, see \cite{Lo}. Another refinement of LERF
for Kleinian groups is GFERF; namely if $G$ is a Kleinian group, then
$G$ is called {\em GFERF} if $G$ is $H$-separable
for every geometrically finite subgroup $H$ of $G$. This has a generalization
when $G$ is a word hyperbolic group; $G$ is called
{\em QCERF} if $G$ is $H$-separable for every finitely generated,
quasi-convex subgroup $H$ of $G$. 
 
We restrict attention to closed orientable hyperbolic $3$-manifolds,
since in this dimension, it is well-known that the fundamental group
of a finite volume, non-compact hyperbolic $3$-manifold or a non-orientable closed
hyperbolic 3-manifold surjects onto $\bf Z$.

\begin{theorem}
\label{main1}
Let $M={\bf H}^3/\Gamma$ be a closed orientable hyperbolic 3-manifold. Assume
that $\Gamma$ has the property that every infinite index, geometrically finite subgroup 
of $\Gamma$ is engulfed in $\Gamma$. 

Then the Lubotzky-Sarnak Conjecture holds for $\Gamma$.\end{theorem}

An immediate corollary of this is (notation as in Theorem \ref{main1}).

\begin{corollary}
\label{lerf}
If $\Gamma$ is LERF, then the Lubotzky-Sarnak Conjecture
holds for $\Gamma$.
\end{corollary}

There is now some
evidence that the fundamental group of any finite volume hyperbolic 3-manifold
is LERF (see \cite{AGM}, \cite{ALR} and \cite{HW} to name a few).
Moreover, Corollary \ref{lerf} was previously
known to hold, if in addition $\Gamma$ contains
a surface subgroup (see \S 3.2 where we give a proof for convenience).

It is also interesting to compare Corollary \ref{lerf}
with the result that if an arithmetic lattice in a semi-simple
Lie group is LERF it cannot have the Congruence Subgroup Property (see
\cite{LR} Chapter 4 for example). It is a consequence of
Clozel's work \cite{Cl} (which
is the culmination of work of many authors) that if an 
arithmetic lattice in a semi-simple
Lie group has the Congruence Subgroup Property it has Property $\tau$.

Another interesting corollary follows from \cite{LaLR} (see  \S 3 for a proof).

\begin{corollary}
\label{arithmetic}
Assume that the fundamental group of every closed hyperbolic 3-manifold is
GFERF. Then, if $\Gamma$ is an arithmetic Kleinian group, $\Gamma$ is
large.
\end{corollary}

It has recently been proved by Agol, Groves and Manning \cite{AGM}
that if every word hyperbolic group is residually finite, then every
word hyperbolic group is QCERF. Combining this with above result, we obtain
the following unexpected conclusion.

\begin{corollary}
\label{resfin}
Assume that every word hyperbolic group is residually finite.
Then every arithmetic Kleinian group is large.
\end{corollary}

Finally we point out that while the Lubotzky-Sarnak
Conjecture remains open, our results have the following consequence
even in the absence of the LERF hypothesis. We let $h(X)$ denote the
Cheeger constant of a Riemannian manifold, possibly with infinite volume. 
When the manifold has finite volume, this is defined to be 
$$h(X) = \inf_S~{  {\Area(S)}\over{\min\{\vol(X_1),\vol(X_2) \}} } $$
where the infimum is taken over all smooth co-dimension one submanifolds
$S$ that separate $X$ into submanifolds $X_1$ and $X_2$. When $X$
has infinite volume, the Cheeger constant is defined to be 
$$h(X) = \inf_S~{  {\Area(S)}\over{\vol(X_1)} } $$
where the infimum is taken over all smooth co-dimension one submanifolds $S$
that bound a compact submanifold $X_1$.

\begin{theorem}
\label{noncompact}
Let $M$ be a closed hyperbolic $3$-manifold. 

Then there is a sequence of (possibly infinite) coverings $M_i$ 
for which $h(M_i) \rightarrow 0$.
\end{theorem}

This result has recently been used by the first author \cite{La2} to
show that nonelementary Kleinian groups which contain a finite
noncyclic subgroup are either virtually free, or contain the
fundamental group of a closed orientable surface of positive genus.  In particular,
co-compact arithmetic Kleinian groups contain surface subgroups.

%
%
%
%
\section{Two preliminary propositions}

Let $N$ be a possibly noncompact complete
Riemannian manifold and $\Delta$ the Laplace-Beltrami operator, 
with sign chosen so that this is a positive operator. Set

$$ \lambda_0(N) = \inf\biggl({{\int_N || \nabla f  ||^2}
                       \over{\int _N f^2}}\biggr),$$
where the infimum is taken over smooth functions $f$ of compact support. 
It is shown in \cite{CY} that $\lambda_0(N)$ is the greatest lower bound of the spectrum 
of $\Delta$ acting on $L^2(N)$. \\[\baselineskip]
\noindent {\bf Remark:}~When $N$ is a closed Riemannian manifold, 
$\lambda_0(N)=0$, and
it is $\lambda_1(N)$ (the first non-zero eigenvalue of $\Delta$) that
is computed by the above infimum, except that $f$ is required
to be orthogonal to the constant functions.

\begin{proposition}
\label{promotingsmallevalue}
Let $M$ be a closed Riemannian manifold, $H$ an infinite
index finitely generated subgroup of $\pi_1(M)$ and $N$ the cover of $M$
corresponding to $H$. Suppose that $\pi_1(M)$ is $H$-separable. 

Then given $\epsilon  > 0$, 
there is a finite sheeted cover $\tilde{M}$ of $M$ for which
$\lambda_1(\tilde{M}) < \lambda_0(N)+\epsilon$.
\end{proposition}

\noindent{\bf Proof:}~ Set $\delta = \epsilon/(1+\epsilon+\lambda_0(N))$.
By \cite{CY}, we may fix some compactly supported function 
$f : N \longrightarrow {\Bbb R}$ for which 

$$ \lambda_0(N) +\delta  > {{\int_N || \nabla f  ||^2}
                       \over{\int _N f^2}}.$$
Choose some compact set $X \subset N$ so that 
${\rm support}(f) \subset {\rm interior}(X)$. 

Since $\pi_1(M)$ is $H$-separable, we may find a finite sheeted
covering, $\widetilde{M}$ of $M$ which is subordinate to $N$ and
for which the compact set $X$ is embedded by the projection 
$ N \longrightarrow \widetilde{M}$ (see \cite{Sc}).  By choosing a larger covering if
necessary, we may arrange that
$${1\over{\vol(\widetilde{M})}}(\int_X f )^2 < \delta \int_X f^2.$$

Define a function $g : \widetilde{M} \longrightarrow {\Bbb R}$ to be
$f$ on $X$ and zero elsewhere. It follows that
$$ \lambda_0(N) +\delta  > 
{{\int_{\widetilde{M}} || \nabla g  ||^2}\over{\int_{\widetilde{M}} g^2}}.$$

We need to adjust the function $g$ slightly, since it is not
orthogonal to the constant functions. This is achieved by replacing 
$g$ by $g^* = g - \alpha$ where $\alpha$ is the constant function 
whose value is $(\int_{\widetilde{M}} g)/\vol(\widetilde{M})$. Then 

$$\int _{\widetilde{M}} (g^*)^2 = \int _{\widetilde{M}} g^2 - 2\int
_{\widetilde{M}} \alpha g + \int _{\widetilde{M}} \alpha^2 = \int
_{\widetilde{M}} g^2 - (\int_{\widetilde{M}} g)^2/{\vol(\widetilde{M})}.$$

Now by construction of $g$, $\int _{\widetilde{M}} g^2 = \int _{X}
f^2$ and $(\int_{\widetilde{M}} g)^2 = (\int_{X} f)^2$, so that the
right hand side of this expression satisfies
$$
\int _{\widetilde{M}} g^2 - (\int_{\widetilde{M}}
g)^2/\vol(\widetilde{M}) = \int _{X} f^2 - (\int_{X}
f)^2/\vol(\widetilde{M}) > (1-\delta) \int_X f^2 = (1-\delta) \int _{\widetilde{M}}
g^2 $$
so that
$(1-\delta)^{-1} >  (\int  _{\widetilde{M}} g^2 )/(\int  _{\widetilde{M}} (g^*)^2 )$. 
\\[\baselineskip]
Now, noting that $\nabla g^* = \nabla g$ we compute

$$
\lambda_1(\widetilde{M}) \leq \int_{\widetilde{M}} || \nabla g^*
||^2 \bigg / \int _{\widetilde{M}} (g^*)^2 = \int_{\widetilde{M}} ||
\nabla g ||^2 \bigg / \int _{\widetilde{M}} (g^*)^2$$
$$ < (\lambda_0(N)+\delta)\int_{\widetilde{M}} g^2  \bigg /   \int _{\widetilde{M}} (g^*)^2 
< (\lambda_0(N)+\delta)/(1-\delta) = \lambda_0(N) + \epsilon$$
as required. \qed

The following result is an analogue of the above proposition, but using
Cheeger constants rather than the first eigenvalue of the Laplacian.

\begin{proposition}
\label{promotingsmallcheeger}
Let $M$ be a closed Riemannian manifold, $H$ an infinite
index finitely generated subgroup of $\pi_1(M)$ and $N$ the cover of $M$
corresponding to $H$. Suppose that $\pi_1(M)$ is $H$-separable. 

Then given $\epsilon  > 0$, 
there is a finite sheeted cover $\tilde{M}$ of $M$ for which
$h(\tilde{M}) < h(N)+\epsilon$.
\end{proposition}

\noindent{\bf Proof:}~ Let $X$ be some compact submanifold of $N$
with zero codimension, and 
such that $\Area(\partial X)/\vol(X) < h(N) + \epsilon$.

Since $\pi_1(M)$ is $H$-separable, we may find a finite sheeted
covering, $\widetilde{M}$ of $M$ which is subordinate to $N$ and
for which the compact set $X$ is embedded by the projection 
$ N \longrightarrow \widetilde{M}$.  By choosing a larger covering if
necessary, we may arrange that $\vol(\widetilde{M}) > 2 \ \vol(X)$.
So,
$$h(\widetilde{M}) \leq \Area(\partial X)/\vol(X) < h(N) + \epsilon$$
as required. \qed

%
%
%
\section{Proof of Theorem \ref{main1}}

In the setting of the fundamental groups of closed Riemannian manifolds 
the definition of Property $\tau$ described in \S 1 is equivalent to the
following (see \cite{Lu} Chapter 4):\\[\baselineskip]
Let $X$ be a closed Riemannian manifold and let
$\Gamma=\pi_1(X)$. Then $\Gamma$ or $X$ has Property $\tau$ if there is
a constant $C>0$ such that $\lambda_1(N)>C$ for all 
finite sheeted covers $N$ of $X$.

\subsection{}
We need the following proposition. (For the definition of Hausdorff
dimension we refer the reader to \cite{Su}.)

\begin{proposition}
\label{fuzzyfree}
Let $M={\bf H}^3/\Gamma$ be a closed hyperbolic 3-manifold. Then $\Gamma$
contains an infinite sequence of finitely generated, free, convex cocompact subgroups
$\{F_j\}$ such that $\lambda_0({\bf H}^3/F_j)\rightarrow 0$.
\end{proposition}
\noindent{\bf Proof:} This is a consequence of results of Sullivan \cite{Su} and L. Bowen \cite{Bo}.

For, it is shown in \cite{Su} that if  $N={\bf H}^3/\Gamma$ is a geometrically finite hyperbolic 
$3$-manifold and $D$ the Hausdorff dimension of the limit set of $\Gamma$, then 
 $\lambda_0(N)=1$ if and only if $D\leq 1$ and otherwise $\lambda_0(N)=D(2-D)$.

Now Bowen shows in \cite{Bo} (actually he shows more than this, but this suffices
for our purpose) that if  $M={\bf H}^3/\Gamma$ is a closed hyperbolic $3$-manifold,
then $\Gamma$ contains an infinite sequence of 
finitely generated, free, convex cocompact subgroups $\{F_j\}$ such that the Hausdorff dimension
of the limit sets of $F_j$ tend to $2$.\qed\\[\baselineskip] 
{\bf Remarks:}\\
{\bf 1.}~By the solution to the Tameness
Conjecture \cite{A} and \cite{CG}, all finitely generated free subgroups
of a cocompact Kleinian group are convex cocompact.
However, Bowen proves that the subgroups $F_j$ are convex co-compact
without appealing to this theorem (see Lemma 5.3 in \cite{Bo}).
\\[\baselineskip]
\noindent{\bf 2.}~Note that it is a consequence of Sullivan's result
above that if $N={\bf H}^3/\Gamma$ is a geometrically finite
hyperbolic 3-manifold and $\Gamma_1$  is a supergroup or 
subgroup of $\Gamma$ of finite index then 
$\lambda_0({\bf H}^3/\Gamma_1)=\lambda_0(N)$.\\[\baselineskip]
We can now complete the proof of our main result.\\[\baselineskip]
\noindent{\bf Proof of Theorem \ref{main1}:}~We begin
with a reduction. We can assume that
all finitely generated subgroups $F$ of $\Gamma$
are geometrically finite. For, if not, then by the solution to the Tameness
Conjecture, $F$ is the fundamental group of a virtual fibre in a fibration over
the circle. It is well known that the Lubotzky-Sarnak
Conjecture holds in this case.

Given Proposition \ref{fuzzyfree}, the remarks following it and 
Proposition \ref{promotingsmallevalue} or \ref{promotingsmallcheeger} it clearly suffices to
prove the following: given a finitely generated, free, convex cocompact subgroup $F$ in $\Gamma$ 
then there is subgroup $F'<\Gamma$ such that $[F':F]<\infty$
and $\Gamma$ is $F'$-separable. This follows immediately from the
engulfing hypothesis using \cite{Lo} 
Theorem 2.7.\qed\\[\baselineskip]
\noindent{\bf Proof of Corollary \ref{arithmetic}:}~This is seen
as follows. Theorem 1.9 of \cite{LaLR} shows that if
every compact 3-manifold with infinite fundamental group does not have
Property $\tau$,  then arithmetic Kleinian groups are large. 

Now assuming the Geometrization Conjecture, this is well known for
compact 3-manifolds which are not hyperbolic (see e.g. the Appendix in
\cite{LaLR}) and
the remaining case is provided by Theorem \ref{main1} (since
GFERF obviously implies the engulfing property for infinite index,
geometrically finite subgroups).\qed\\[\baselineskip]
Finally, Theorem \ref{noncompact} follows from the non-compact version
of Cheeger's inequality $\lambda_0(N)\geq h(N)^2/4$ applied to
Proposition \ref{fuzzyfree}.\qed

\subsection{} We include the following argument for convenience, 
and to emphasize Corollary \ref{arithmetic}.

\begin{theorem}
\label{lerftolargeifasurface}
Let $M={\bf H}^3/\Gamma$ be a closed hyperbolic 3-manifold, and assume
that $\Gamma$ contains the fundamental group of a closed surface of genus
at least $2$. Then if $\Gamma$ is LERF, then
$\Gamma$ is large.\end{theorem}

\noindent{\bf Proof:}~The surface subgroup corresponds to a closed 
incompressible surface immersed in $M$. If $S$ is geometrically infinite,
then it is a virtual fiber in a fibration over the circle. By passing to
a finite sheeted cover, it follows from \cite{CLR} that $\Gamma$ must
also contain a closed quasi-Fuchsian surface subgroup. 
Thus we now work with $F$ a 
quasi-Fuchsian surface subgroup of $\Gamma$. 

Using the LERF assumption, we invoke Scott's result \cite{Sc} to pass
to a finite sheeted cover $M_1={\bf H}^3/\Gamma_1$ so that $M_1$
contains a closed embedded quasi-Fuchsian surface with covering group
$F$. This determines a free product with amalgamation decomposition 
$A *_F B$ or HNN-extension $A*_F$ for $\Gamma_1$. One can now arrange a
surjective homomorphism onto an amalgam of finite groups to finish the
proof (see for example \cite{Lu1}).  Note that, since $F$ is quasi-Fuchsian,
the amalgam cannot be ${\bf Z}/2{\bf Z} * {\bf Z}/2{\bf Z}$.\qed

\section{Final Comments}

\noindent{\bf 1.}~It is important in Bowen's proof that there exist
discrete, convex compact, free Kleinian groups, where the Hausdorff
dimension of their limit set is arbitrarily close to $2$. That there
are examples of purely
hyperbolic free subgroups whose limit set is the entire sphere at
infinity seems to have first been established by Greenberg \cite{Gr}
(as points on the boundary of Schottky space which
are limits of convex cocompact groups). That the
Hausdorff dimension of the limit sets of these convex cocompact groups 
get arbitrarily close to $2$
can be seen from Corollary 7.8 of \cite{Mc} for example.

The existence of analogous subgroups in
$\SO(n,1)$ for $n\geq 4$ is as yet unknown. 
Their existence, together with the known
generalization of Sullivan's result \cite{Su} to higher dimensions
would prove that LERF
implies the Lubotzky-Sarnak Conjecture for higher dimensional hyperbolic
manifolds.\\[\baselineskip]
\noindent{\bf 2.}~It is interesting to contrast Proposition \ref{fuzzyfree}
with what happens, for instance for (free) subgroups of cocompact lattices
in $\Sp(n,1)$, $n\geq 2$. If $\Gamma$ is such a lattice then
it has Property T. It is shown in \cite{Br} (see Theorem 3), in 
contrast to Proposition \ref{fuzzyfree},
that if $\Delta$ is a subgroup of $\Gamma$, which is either
finite or infinite index, there is a spectral gap for the smallest
non-zero eigenvalue of the Laplacian. 

A similar result was established in \cite{Co} for the Hausdorff
dimension; namely that any infinite index
convex cocompact subgroup of $\Gamma$ (as above) has Hausdorff co-dimension
of its limit set being at least $2$.

Neither LERF nor the Congruence Subgroup Property are known for any example
in this setting.\\[\baselineskip]
\noindent{\bf 3.}~Other situations where LERF is used to imply large
were recently given in \cite{Bu}.

%
%
%
%

\noindent Mathematical Institute,\\ University of Oxford,\\
24-29 St Giles',\\ Oxford OX1 3LB, UK.

\noindent Email:lackenby@maths.ox.ac.uk\\[\baselineskip]
 Department of Mathematics,\\ University of California,\\ Santa Barbara, CA
93106, USA.

\noindent Email:~long@math.ucsb.edu\\[\baselineskip]
 Department of Mathematics,\\
 University of Texas,\\
 Austin, TX 78712, USA.

\noindent Email:~areid@math.utexas.edu\\

\end{document}